\documentclass{article}

\usepackage{amsmath,amsfonts,a4wide,latexsym,theorem,amssymb}

{\theorembodyfont{\slshape}
\newtheorem{proposition}{Proposition}[section]
\newtheorem{lemma}[proposition]{Lemma}
\newtheorem{corollary}[proposition]{Corollary}
\newtheorem{theorem}[proposition]{Theorem}}

{\theorembodyfont{\upshape}
\newtheorem{remark}[proposition]{Remark}}

\newcommand\qed{\hfill$\square$}

\newcommand\C{{\mathbb C}}

\newcommand\Q{{\mathbb Q}}
\newcommand\Z{{\mathbb Z}}

\newcommand\tilf{{\tilde f}}
\newcommand\tilP{{\widetilde P}}

\newcommand\tilrho{{\tilde \rho}}
\newcommand\tilvarrho{{\tilde \varrho}}

\newcommand\bfa{{\mathbf a}}
\newcommand\bfj{{\mathbf j}}

\newcommand\calA{{\mathcal A}}
\newcommand\calD{{\mathcal D}}

\newcommand\norm{{\mathcal N}}

\newcommand\vheight{{\mathrm h}}

\newcommand\aheight{{\mathrm{h_a}}}
\newcommand\pheight{{\mathrm{h_p}}}

\newcommand\ord{{\mathrm {ord}}}

\newcommand\Mu{{\mathrm M}}

\makeatletter

\renewcommand*\l@section[2]{%
  \ifnum \c@tocdepth >\z@
    \addpenalty\@secpenalty
    \addvspace{0.3em \@plus\p@}%
    \setlength\@tempdima{1.5em}%
    \begingroup
      \parindent \z@ \rightskip \@pnumwidth
      \parfillskip -\@pnumwidth
      \leavevmode \bfseries
      \advance\leftskip\@tempdima
      \hskip -\leftskip
      #1\nobreak\hfil \nobreak\hb@xt@\@pnumwidth{\hss #2}\par
    \endgroup
  \fi}
  
 \makeatother
 
\title{Remarks on Eisenstein}

\author{Yuri Bilu (Bordeaux), Alexander Borichev (Marseille)}

\date{December 21, 2011}

1

\begin{document}

\maketitle

\begin{abstract}
We obtain a fully explicit quantitative version of the Eisenstein theorem on algebraic power series which is more suitable for certain applications than the existing version due to Dwork, Robba, Schmidt and van der Poorten. We also treat ramified series and Laurent series, and we demonstrate some applications; for instance, we estimate the discriminant of the number field generated by the coefficients. 
\end{abstract}

{\footnotesize

\tableofcontents

\bigskip

\begin{flushright}
\textit{Dedicated to the memory of Alf van der Poorten}
\end{flushright}

}

\section{Introduction}
Let ${f(z)=a_0+a_1z+a_2z^2+\ldots\in \bar\Q[[z]]}$ be an algebraic power series. ``Algebraic'' means here that~$f$ is algebraic over the field of rational functions $\bar\Q(z)$. It is well-known (and easy to see) that the coefficients of an algebraic power series belong to a finite extension of~$\Q$.  The classical theorem attributed to Eisenstein asserts that there exists  a non-zero rational integer~$T$ such that  ${T^{k+1}a_k}$ is an algebraic integer for all~$k$. 

The fact that $f(z)$ is algebraic means that there exists a polynomial ${P(z,w)\in\bar\Q[z,w]}$ such that ${P\bigl(z,f(z)\bigr)=0}$. 
Eisenstein himself made his observation in~\cite{Ei52} assuming the non-vanishing condition ${P'_w(0,0)\ne0}$ (the ``implicit function theorem'' condition). In this case the proof is simpler and one can have a good estimate for the ``Eisenstein constant''~$T$. The first proof in the general case was, probably, due to Heine \cite[pages 50--53]{He78}, who observed that for sufficiently large~$N$ the ``tail series'' ${\sum_{k=N}^\infty a_kz^{N-k}}$ satisfies a polynomial equation with non-vanishing $w$-derivative at the origin. 

Assume that the coefficients of the series $f(z)$ belong to a number field~$K$. Denote by~$M_K$ the set of absolute values  of~$K$ (normalized to extend the standard  absolute values of~$\Q$).  The Eisenstein theorem is essentially equivalent to the following statement: for every ${v\in M_K}$ the $v$-adic norms $|a_k|_v$ grow at most exponentially in~$k$, and for all but finitely many~$v$ we have  ${|a_k|_v\le 1}$ for all~$k$. (Strictly speaking, the classical Eisenstein theorem refers only to finite absolute values, but it is natural to include infinite absolute values as well.)

It will be convenient to use the notion of $M_K$-divisor. By an \textsl{$M_K$-divisor} we mean associating to every ${v\in M_K}$ a positive real number~$A_v$ with the following property: for all but finitely many~$v$ we have ${A_v=1}$. 
Then Eisenstein theorem simply means that there exists an $M_K$-divisor such that ${|a_k|_v\le A_v^{k+1}}$ for every ${v\in M_K}$.

Eisenstein theorem has many applications, most notably in Diophantine analysis \cite{Ab08,BC70,Bi97,BSS12,Co70,CZ03,Fu03,HS83,LP04,SS08,Sc92}, to mention only some of them. 
For these  applications  one often needs a quantitative version of the theorem  in terms of the polynomial $P(z,w)$ above. For instance, in the Diophantine analysis one often needs to construct explicitly rational functions on algebraic curves with prescribed poles and zeros, see \cite{Co70,Sc91}. The Eisenstein theorem in an indispensable tool for such constructions, and to make them explicit one needs quantitative version of the qualitative statement above. First such versions were given by Coates~\cite{Co70} and Hilliker-Strauss~\cite{HS83}, who used the ``tail series'' trick described above.

Schmidt~\cite{Sc90} suggested a much more efficient approach: to estimate the $v$-adic norms of the coefficients  in terms of the $v$-adic convergence radius. Using the previous work of Dwork and Robba~\cite{DR79} Schmidt obtained optimal estimates when the underlying prime~$p_v$ satisfies ${n< p_v\le \infty}$, where ${n=\deg_wP(z,w)}$. Schmidt conjectured that the same was possible for ${p_v\le n}$ as well (his own estimates in this case were not optimal).  This was confirmed by Dwork and van der Poorten~\cite{DP92}.

The following theorem is a compilation of the results from \cite{DR79,DP92,Sc90}, maid in \cite[Theorem~2.2]{Bi97}.
We define the \textsl{height} of an $M_K$-divisor ${\calA=(A_v)_{v\in M_K}}$ by 
$$
\vheight(\calA)= d^{-1}\sum_{v\in M_K} d_v\log^+A_v,\qquad \log^+=\max\{\log, 0\}. 
$$
where ${d=[K:\Q]}$ is the absolute degree of~$K$ and ${d_v=[K_v:\Q_v]}$ is the absolute local degree of~$v$. 

\begin{theorem}[Dwork, Robba, Schmidt, van der Poorten]
\label{tdrsvdp}
Let~$K$ be a number field and ${P(z,w)\in K[z,w]}$ an irreducible polynomial of degrees 
${\deg_zP=m}$ and  ${\deg_wP=n}$. Further, let
${f(z)=\sum_{k=0}^\infty a_kz^k \in \bar K[[z]]}$ be 
a power series  satisfying ${P\bigl(z,f(z)\bigr)=0}$. Then there exist $M_K$-divisors ${\calA'=(A_v')_{v\in M_K}}$ and ${\calA=(A_v)_{v\in M_K}}$ such that 
$$
|a_k|_v\le A_v'A_v^{m+k}\qquad (k=0,1,2,\ldots)
$$ 
for any ${v\in M_K}$  anyhow extended to $\bar K$,  and such that 
\begin{equation}
\label{edrsv}
\vheight(\calA') \le \pheight(P)+O(\log n),  \quad  \vheight(\calA)\le (2n-1)\pheight(P) + O\bigl(n\log(2mn)\bigr). 
\end{equation}
\end{theorem}

Here $\pheight(P)$ is the usual projective height of the polynomial~$P$, see Section~\ref{sprelim}.

As Dwork and van der Poorten explain in the introduction to their article~\cite{DP92}, one cannot improve on the term ${(2n-1)\pheight(P)}$. So, the bound for $\vheight(\calA)$ in~\ref{edrsv}  is best possible. Still, there is one unsatisfactory point in Theorem~\ref{tdrsvdp}: the term~$m$ in the estimate ${|a_k|_v\le A_v'A_v^{m+k}}$. This term originates from Schmidt's Lemma~2 in~\cite{Sc90}. It is implicit in Schmidt's article and indicated in~\cite{Bi97} that~$m$ can be replaced by $\ord_zp_n(z)$, where $p_n(z)$ is the leading coefficient of $P(z,w)$ viewed as polynomial in~$w$ over $K[z]$. Still, one expects, and for certain applications one needs a bound ${|a_k|_v\le A_v'A_v^{k}}$ with divisors~$\calA$ and~$\calA$ satisfying~\eqref{edrsv} or similar. 

The principal purpose of this article is obtaining such a bound. The price we have to pay here is a slightly weaker estimate for the height of~$\calA$, with ${2n-1}$ replaced by ${3n-1}$. The following theorem is proved in Section~\ref{sglo}. 

\begin{theorem}
\label{tmain}
In the set-up of Theorem~\ref{tdrsvdp} there exist effective $M_K$-divisors ${\calA'=(A_v')_{v\in M_K}}$ and ${\calA=(A_v)_{v\in M_K}}$ such that 
$$
|a_k|_v\le A_v'A_v^k\qquad (k=0,1,2,\ldots)
$$ 
for any ${v\in M_K}$  anyhow extended to $\bar K$,  and such that 
\begin{equation}
\label{eeisenreg0}
\vheight(\calA') \le \pheight(P)+\log 3,  \quad  \vheight(\calA)\le (3n-1)\pheight(P) + 3n\log (mn) +7n. 
\end{equation}
\end{theorem}

We also obtain more general results, on  ramified series and Laurent series (they cannot be reduced to Theorem~\ref{tmain} by a simple variable change).

\paragraph{Plan of the article}
In Section~\ref{sprelim} we introduce some (mainly standard) terminology. In Section~\ref{sestim} we collect various auxiliary results to be used in the main body of the article.

Section~\ref{sloreg} is the technical heart of the article. We obtain therein a required version of the already mentioned Schmidt's Lemma~2 from~\cite{Sc90}, and deduce a $v$-adic upper bound for the coefficients of the series. In Section~\ref{slogen} we extend the results of the previous section to general ramified Laurent series.

Theorem~\ref{tmain}, as well as its ramified Laurent generalization are deduced in Section~\ref{sglo} more or less straightforwardly from the results of Sections~\ref{sloreg} and~\ref{slogen}. Some immediate application are given as well.

In Section~\ref{sa0} we obtain a slightly different type of Eisenstein theorem, where the estimate involves the initial coefficient~$a_0$. Finally, in Section~\ref{sfield} we give an explicit upper bound for the discriminant of the number field generated by the coefficients.

\paragraph{Acknowledgments}
We thank Michel Matignon for suggesting an elegant proof of Lemma~\ref{lbc}, and Andrea Surroca for valuable discussions. 

Yuri Bilu was supported by the ANR project HAMOT, and by the Swiss National Foundation \textit{Ambizione} fund PZ00P2\_121962. He thanks the University of Basel and Andrea Surroca for hospitality in November-December 2011, when a substantial part of this work was done.

\section{Notation and Conventions}
\label{sprelim}

In this section we make general definitions and state some conventions assumed throughout the article. 

Le~ $P$ be a polynomial over some field. For any non-zero element~$\alpha$ of the field, the polynomial $\alpha P$ will be called a \textsl{normalization} of~$P$. A normalization will be called \textsl{moderate} if one of its coefficients is~$1$.

Let~$P$ be a polynomial over a field supplied with an absolute value ${|\cdot|}$. We denote by $|P|$ the maximal absolute value of the coefficients of~$P$. 

We denote by $\bar K$ the algebraic closure of the field~$K$. 

Let~$K$ be a number field of degree ${d=[K:\Q]}$. We denote by~$M_K$ the set of  absolute values of~$K$ normalized to extend the standard absolute values of~$\Q$. That is, for ${v\in M_K}$ we have ${|p|_v=p^{-1}}$ if $v\mid p<\infty$, and ${|2012|_v=2012}$ if ${v\mid\infty}$. With this convention the product formula reads 
$$
\prod_{v\in M_K} |\alpha|_v^{d_v}=1 \qquad (\alpha\in K^\ast),
$$ 
where ${d_v=[K_v:\Q_v]}$ is the the local degree of~$v$.

Given a vector ${\bfa=(a_1, \ldots, a_m)}$ with algebraic entries, we define its \textsl{projective height} and its \textsl{affine height} as
$$
\pheight(\bfa)= d^{-1}\sum_{v\in M_K}d_v\log\max_{1\le k \le m}|a_k|_v, \qquad \aheight(\bfa)= d^{-1}\sum_{v\in M_K}d_v\log^+\max_{1\le k \le m}|a_k|_v,
$$
where ${\log^+=\max\{\log,0\}}$ and~$K$ is a number field containing the entries of~$\bfa$, the degrees~$d$ and~$d_v$ being as above. It is well-known (and very easy to see) that  both $\pheight(\bfa)$ and $\aheight(\bfa)$ depend only on the vector~$\bfa$, but not on the particular choice of the field~$K$.

Given a polynomial~$P$ with algebraic coefficients, we define its heights $\pheight(P)$ and $\aheight(P)$ as the corresponding heights of the vector of its non-zero coefficients. If the polynomial~$P$ has~$1$ as one of its coefficients (that is,~$P$ is \textsl{moderately normalized}  in the terminology introduced above), then ${\pheight(P)=\aheight(P)}$. 

Finally, in this article the letter~$e$ is used exclusively for ramification indices; it never denotes the constant ${2.718\dots}$

\section{Some Estimates}
\label{sestim}

In this section we collect  some very simple estimates which will be used in the article.

\subsection{Local Estimates}

In this subsection~$K$ is an algebraically closed field supplied with an absolute value~${|\cdot|}$. For a polynomial~$P$ over~$K$ we denote by $|P|$ the maximum of absolute values of the coefficients of~$P$.

\begin{proposition}
\label{pzw}
\begin{enumerate}
\item
\label{iw}
Let $P(w_1, \ldots, w_\ell)$ be a polynomial over~$K$.   
Put ${n_i=\deg_{w_i}P}$. Then at any point ${(w_1, \ldots, w_\ell)\in K^\ell}$  we have
$$
\left|P(w_1, \ldots, w_\ell)\right|\le A|P|\prod_{i=1}^\ell\max\{1,|w_i|\}^{n_i},
$$
where
${A=(n_1+1)\cdots(n_\ell+1)}$ in the archimedean case, and ${A=1}$ 
in the non-archimedean case.

\item
\label{izw}
Let $P(z_1,\ldots, z_k, w_1, \ldots, w_\ell)$ be a polynomial over~$K$.   
Put ${n_i=\deg_{w_i}P}$. Then at any point ${(z_1,\ldots, z_k, w_1, \ldots, w_\ell)\in K^{k+\ell}}$ such that ${|z_1|, \ldots |z_k|<1}$ and ${|w_1|, \ldots, |w_\ell|>1}$ we have
$$
\left|P(z_1,\ldots, z_k, w_1, \ldots, w_\ell)\right|\le 
\begin{cases}
\displaystyle|P|\prod_{i=1}^k\frac1{1-|z_i|}\prod_{i=1}^\ell\frac{|w_i|^{n_i}}{1-|w_i|^{-1}}, &\text{$|\cdot|$ is archimedean},\\
\displaystyle|P|\prod_{i=1}^\ell|w_i|^{n_i}, &\text{$|\cdot|$ is non-archimedean}.
\end{cases}
$$
\end{enumerate}
\end{proposition}

\paragraph{Proof}
Item~\ref{iw} and the non-archimedean case of item~\ref{izw} are obvious. In the archimedean case of item~\ref{izw}, denoting ${m_i=\deg_{z_i}P}$, we have
{\small%
$$
\left|P(z_1,\ldots, z_k, w_1, \ldots, w_\ell)\right| \le |P| \prod_{i=1}^\ell|w_i|^{n_i} \prod_{i=1}^k\left(1+|z_i|+ \cdots+|z_i|^{m_i}\right)\prod_{i=1}^\ell\left(1+|w_i|^{-1}+ \cdots+|w_i|^{-n_i}\right).
$$}%
Replacing finite sums by infinite sums, the result follows.
\qed

\begin{proposition}
\label{pboundroot}
Let~$\alpha$ be a non-zero root of ${P(z)=a_nz^n+\cdots+a_0\in K[z]}$ (with ${a_n\ne 0}$) and let~$m$ be the order of~$0$ as a root of~$P$. Then 
\begin{align*}
\displaystyle\frac{|a_m|}{2|P|}\le |\alpha|\le\displaystyle\frac{2|P|}{|a_n|}, \qquad&\text{$|\cdot|$ is archimedean},\\
\displaystyle\frac{|a_m|}{|P|}\le |\alpha|\le \displaystyle\frac{|P|}{|a_n|}, \qquad&\text{$|\cdot|$ is non-archimedean}.
\end{align*}
\end{proposition}

\paragraph{Proof}  For the upper estimates we may assume that ${|\alpha|\ge2}$ in the archimedean case and ${|\alpha|>1}$ in the non-archimedean case. Using item~\ref{izw} of Proposition~\ref{pzw} we obtain
$$
 |\alpha|^n = \left|\frac{a_{n-1}}{a_n}\alpha^{n-1}+ \cdots+\frac{a_0}{a_n}\right| \le
\begin{cases}
\displaystyle\frac{|P|}{|a_n|}\frac{|\alpha|^{n-1}}{1-|\alpha|} \le \frac{2|P|}{|a_n|}|\alpha|^{n-1}, &\text{$|\cdot|$ is archimedean},\\
\displaystyle\frac{|P|}{|a_n|}|\alpha|^{n-1}, &\text{$|\cdot|$ is non-archimedean},
\end{cases}
$$
which proves the upper estimates. The lower estimates follow upon replacing~$\alpha$ by~$\alpha^{-1}$. \qed

\medskip

The following proposition will be used only in the case ${k=\ell=1}$, but we prefer to state it in the general case for the sake of further applications.

\begin{proposition}
\label{ptranslate}
\begin{enumerate}
\item
\label{ideriv}
Let $P(w_1, \ldots, w_\ell)$ be a polynomial over~$K$.   
Put ${n_i=\deg_{w_i}P}$. For a multy-index ${\bfj=(j_1, \ldots, j_\ell)}$ denote by $D_\bfj$ the differential operator 
$$
D_\bfj= \frac1{j_1!\cdots j_\ell!}\frac{\partial^{j_1+\cdots+j_\ell}}{\partial w_1^{j_1}\cdots \partial w_\ell^{j_\ell}}
$$ 
Then
${\bigl|D_\bfj P\bigr| \le 2^{n_1+\cdots+n_\ell}|P|}$ in the archimedeal case, and  ${\bigl|D_\bfj P\bigr| \le |P|}$ in the non-archimedean case. 

\item
\label{itranslate}
Let $P(z_1,\ldots, z_k, w_1, \ldots, w_\ell)$ be a polynomial over~$K$, and ${\alpha_1, \ldots, \alpha_\ell \in K}$. Put ${n_i=\deg_{w_i}P}$. Then the polynomial ${Q(z_1,\ldots, z_k, w_1, \ldots, w_\ell)=P(z_1,\ldots, z_k, w_1+\alpha_1, \ldots, w_\ell+\alpha_\ell)}$ 
satisfies
$$
|Q|\le A|P|\prod_{i=1}^\ell\max\{1,|\alpha_i|\}^{n_i},
$$
where
${A=(n_1+1)\cdots(n_\ell+1)2^{n_1+\cdots+n_\ell}}$ in the archimedean case, and ${A=1}$ 
in the non-archimedean case.
\end{enumerate}
\end{proposition}

\paragraph{Proof}
In item~\ref{ideriv} each coefficient of the polynomial $D_\bfj P$ is equal to a coefficient of~$P$ multiplied by a product of binomial coefficients ${\binom {\nu_1}{j_1}\cdots \binom{\nu_\ell}{j_\ell}}$, where ${\nu_i\le n_i}$. The absolute value of this product does not exceed $2^{n_1+\cdots+n_\ell}$ in the archimedean case, and~$1$ in the non-archimedean case. This proves item~\ref{ideriv}. 

In item~\ref{itranslate} every coefficients of~$Q$ is of the form ${D_\bfj(0, \ldots, 0, \alpha_1, \ldots, \alpha_\ell)}$. Hence item~\ref{itranslate} is a consequence of item~\ref{ideriv} of this proposition  and item~\ref{iw} of Proposition~\ref{pzw}.

\subsection{Global Estimates}
\label{ssglob}

\begin{proposition}
\label{ptrah}
Let $P(z_1,\ldots, z_k, w_1, \ldots, w_\ell)$ be a polynomial with algebraic coefficients, and ${\alpha_1, \ldots, \alpha_\ell}$ are algebraic numbers. Put ${n_i=\deg_{w_i}P}$. Then the polynomial 
$$
Q(z_1,\ldots, z_k, w_1, \ldots, w_\ell)=P(z_1,\ldots, z_k, w_1+\alpha_1, \ldots, w_\ell+\alpha_\ell)
$$
satisfies
$$
\pheight(Q) \le \pheight(P) + \sum_{i=1}^\ell \bigl(n_i\aheight(\alpha) +n_i\log 2+\log (n_i+1)\bigr).  
$$
\end{proposition}

\paragraph{Proof} The coefficients of~$P$ and the numbers ${\alpha_1, \ldots, \alpha_\ell}$ belong to some number field~$K$. Then the desired statement is an  immediate  consequence of item~\ref{itranslate} of Proposition~\ref{ptranslate}. \qed

\medskip

Given a polynomial $P(z,w)$, we denote by $R_P(z)$ the $w$-resultant of~$P$ and~$P'_w$. 
\begin{proposition}
\label{prp}
Let ${P(z,w)}$ be a polynomial with algebraic coefficients of $z$-degree~$m$ and $w$-degree~$n$. Then
$$
\pheight(R_P)\le (2n-1)\pheight(P)+(2n-1)\log\left((m+1)(n+1)\sqrt{n}\right).
$$
\end{proposition}
For the proof see Schmidt \cite[Lemma~4]{Sc90}. 

\begin{proposition}
\label{pmahler}
Let ${P(z)}$ be a polynomial with algebraic coefficients of degree~$m$. 

\begin{enumerate}
\item
\label{ioneroot}
Let~$\alpha$ be a root of~$P$. Then ${\aheight(\alpha)\le \pheight(P)+ \log 2}$. 

\item
\label{iallroots}
Let ${\alpha_1, \ldots, \alpha_m}$ be the roots of~$P$ counted with multiplicities. Then 
$$
\aheight(\alpha_1)+\cdots +\aheight(\alpha_m)\le \pheight(P)+\log (m+1).
$$
\end{enumerate}
\end{proposition}

\paragraph{Proof}
Item~\ref{ioneroot} is an immediate consequence of Proposition~\ref{pboundroot}. Item~\ref{iallroots} a classical result of Mahler, see, for instance, \cite[Lemma 3]{Sc90}.\qed

\section{Local Eisenstein Theorem: the Regular Case}
\label{sloreg}

In this section~$p$ is a prime number or ${p=\infty}$, and~$\C_p$ is the completion of the algebraic closure of~$\Q_p$; in particular, ${\C_\infty=\C}$. Since~$p$ is fixed, we may write $|\cdot|$ instead of $|\cdot|_p$.

Let~$\rho$ be a positive real number and let
\begin{equation}
\label{eserf}
f(z)=a_0+a_1z+a_2z^2+\ldots \in \C_p[[z]]
\end{equation}
be a power series converging in the circle ${|z|<\rho}$ and  satisfying the  equation 
${P\bigl(z,f(z)\bigr)=0}$,  where ${P(z,w)\in \C_p[z,w]}$
is a non-zero polynomial. The purpose of this section is to estimate the coefficients~$a_k$ in terms of the polynomial~$P$. First of all let us estimate $f(z)$ in some circle contained in the circle of convergence. 
The principal novelty of our estimate, as compared with that of Schmidt \cite[Lemma~2]{Sc90} is that it does not depend on the $z$-degree of $P(z,w)$.

Recall that we denote by $|P|$ the maximum of absolute values of the coefficients of~$P$.

\begin{proposition}
\label{pbor}
In the above set-up, assume that the polynomial $P(z,w)$ is not divisible by~$z$, so that the polynomial ${q(w):=P(0,w)}$ is non-zero. Assume further that the polynomial $q(w)$ is monic:
$$
q(w)=w^{\deg q}+ \text{terms of lower degree}.
$$ 
Put ${n=\deg_wP}$ and 
$$
\varrho =
\begin{cases}
\min\left\{\rho, (6|P|)^{-n}\right\}, & p=\infty,\\
\min\left\{\rho, |P|^{-n}\right\}, & p<\infty.
\end{cases}
$$
Then for ${|z| < \varrho}$ we have ${|f(z)| \le 3|P|}$ if ${p=\infty}$, and ${|f(z)| \le |P|}$ if ${p<\infty}$.
\end{proposition}

The non-archimedean part of the proof requires a lemma, that may be viewed as a $p$-adic analogue of the ``Bolzano-Cauchy theorem''  from the elementary analysis. Recall that for ${p<\infty }$ the possible absolute values of non-zero elements of~$\C_p$ are exactly the rational powers of~$p$:
$$
\{|z| : z\in \C_p^\ast\}= p^\Q,
$$
where ${p^\Q=\{p^a: a\in \Q\}}$. 
\begin{lemma}
\label{lbc}
Assume that ${p<\infty}$, and let~$f$ be defined by the series~\eqref{eserf} convergent in the circle ${|z|<\rho}$. Let ${s_1, s_2\in p^\Q}$ be such that ${s_1<s_2}$  and  there exist ${z_1,z_2}$ in the circle  with ${|f(z_i)|=s_i}$. Then for any  ${s\in [s_1,s_2]\cap p^\Q}$ there exists~$z$ in the same circle such that ${|f(z)|=s}$. 
\end{lemma}

\paragraph{Proof}
Shifting the variable we may assume that ${z_1=0}$. Convergence of the series~\eqref{eserf} implies that for a  non-negative ${r<\rho}$ we have 
${\lim_{k\to\infty}|a_k|r^k =0}$.
Hence the quantity 
$$
\Mu(r)= \max \bigl\{|a_k|r^k : k=0, 1 , \ldots\bigr\}
$$
is well-defined for  ${r\in [0, \rho)}$. 
Let us show that for ${r\in [0,\rho)\cap p^\Q}$ 
\begin{equation}
\label{emax}
\max\bigl\{|f(z)| : |z|\le r\bigr\}= \Mu(r).
\end{equation}
Denote by~$\ell$ the biggest~$k$ with ${|a_k|r^k=\Mu(r)}$,  and by $f_\ell(z)$ the $\ell$-th partial sum of the series~\eqref{eserf}:
$$
f_\ell(z)=a_0+a_1z+\cdots+a_\ell z^\ell.
$$ 
The ``$\le$''-inequality in~\eqref{emax} is obvious, and to prove the opposite inequality, we must find~$z$ with ${|z|=r}$  for which ${|f_\ell(z)|=\Mu(r)}$. Putting ${g(z)=\Mu(r)^{-1}f_\ell(rz)}$, we  must find~$z$ with ${|z|=1}$  for which ${|g(z)|=1}$. By the definition of $M(r)$, the polynomial $g(z)$ has coefficients in the local ring of~$\C_p$, and its reduction $\bar g(z)$ modulo the maximal ideal is a non-zero polynomial over the residue field of~$\C_p$.  We are left with the following task: find a non-zero element of the residue field which is not a root of $\bar g(z)$.  And this is always possible, because the residue field  is infinite. This proves~\eqref{emax}.

{\sloppy

Since $\Mu(r)$ is a continuous real function on the interval $[0,\rho)$, satisfying  ${\Mu(0)=s_1}$ and ${\Mu(|z_2|)\ge s_2}$, for any ${s\in [s_1,s_2]}$ there exists ${r\in [0,\rho)}$ with ${\Mu(r)=s}$. If~$s$ is a rational power of~$p$, then so is~$r$, and~\eqref{emax} implies that ${f(z)=s}$ for some $z$ in the circle. \qed

}
\medskip

(We owe this elegant argument to Michel Matignon.)

\paragraph{Proof of Proposition~\ref{pbor}}
To simplify the notation, put ${A=|P|}$. It might be worth noticing that ${A\ge1}$ (because $P(z,w)$ has~$1$ as one of the coefficients) and that ${\varrho\le 1}$. Both inequalities will be repeatedly use in the proof. 

The  argument splits into two cases, ${p=\infty}$ and ${p<\infty}$, the proofs being similar, but not identical. 

Assume first that ${p=\infty}$.
Then ${|f(0)|\le 2A}$ by Proposition~\ref{pboundroot}, because $f(0)$ is a root of the monic polynomial $q(w)$.
Hence, if our statement is not true, then there exists~$z$ in the circle ${|z|< \varrho}$ such that ${3A <|f(z)|\le 4A}$. Writing
\begin{equation}
\label{equs}
P(z,w)=q(w)+zQ(z,w)
\end{equation}
and using item~\ref{izw} of Proposition~\ref{pzw}, we obtain for such~$z$ the estimates 
$$
\left|q\bigl(f(z)\bigr)\right|\ge |f(z)|^{\deg q} - A \frac{|f(z)|^{\deg q-1}}{1-|f(z)|^{-1}}\ge
|f(z)|^{\deg q} - \frac32A |f(z)|^{\deg q-1}
>\frac32 A
$$
and
$$
\left|Q\bigl(z,f(z)\bigr)\right|\le A\frac1{1-|z|}\frac{|f(z)|^n}{1-|f(z)|^{-1}}< 2A (4A)^n.
$$
Hence
$$
\left|P\bigl(z,f(z)\bigr)\right|> \frac32 A - 2A (4A)^n(6A)^{-n} >0,
$$
a contradiction. This proofs the statement for ${p=\infty}$.

{\sloppy

In the case ${p<\infty}$ Proposition~\ref{pboundroot} gives ${|f(0)|\le A}$.  If our statement is not true, then there  exists~$z'$ with ${|z'|<\varrho}$ such that ${|f(z')|>A}$. Pick~$r$ such that ${|z'|<r<\varrho}$; notice that ${r<A^{-n}}$. Lemma~\ref{lbc} (with~$r$ instead of~$\rho$) implies the existence of~$z$ with ${|z|<r}$ and ${A<\bigl|f(z)\bigr|<(rA)^{-1/(n-1)}}$. Straightforward estimates give us 
$$
\left|q\bigl(f(z)\bigr)\right|\ge |f(z)|> \bigl|rAf(z)^n\bigr| > \left|zQ\bigl(z,f(z)\bigr)\right|,
$$
whence ${P\bigl(z,f(z)\bigr)\ne0}$, a contradiction. This completes the proof for ${p<\infty}$ as well. 
\qed

}

\begin{corollary}
\label{ccoef}
In the set-up of Proposition~\ref{pbor}, the coefficients of the series~\eqref{eserf} satisfy
$$
|a_k|\le 
\begin{cases}
3|P|\varrho^{-k}, & p=\infty, \\
|P|\varrho^{-k}, & p<\infty.
\end{cases}
$$
\end{corollary}

\paragraph{Proof}
It suffices to show that for any positive ${r<\rho}$ we have
\begin{equation}
\label{eakle}
|a_k| \le r^{-k}\sup\bigl\{|f(z)|: |z|=r\bigr\}. 
\end{equation}
In the case ${p=\infty}$ this follows from
$$
a_k=\frac1{2\pi i}\int_{|z|=r}f(z)z^{-k-1}dz.
$$
In the case ${p<\infty}$ this follows from~\eqref{emax}. \qed

\medskip

To make all this work, we need a lower estimate for the convergence radius~$\rho$ in terms of the polynomial~$P$. This is the principal contents of the work of Dwork, Robba, Schmidt and van der Poorten \cite{DR79,DP92,Sc90}. 
Given a polynomial ${F(z)\in \C_p(z)}$, we denote by $\sigma(F)$  the smallest absolute value of a non-zero root of~$F(z)$:
\begin{equation}
\label{esigma}
\sigma(F)=\min\{|\alpha|: F(\alpha)= 0,  \ \alpha\ne0\}.
\end{equation}
Call a polynomial $P(z,w)$ \ \textsl{$w$-separable} if 
it is not divisible by a square of a polynomial of positive $w$-degree, and denote by $R_P(z)$  the $w$-resultant of $P(z,w)$ and $P'_w(z,w)$; this latter is a non-zero polynomial if $P(z,w)$ is $w$-separable.
\begin{theorem}[Dwork, Robba, Schmidt, van der Poorten]
\label{tdrsv}
Assume that~$P$ is $w$-separable. Then the series~\eqref{eserf} converges for ${|z|<\sigma(R_P)}$ if ${n<p\le \infty}$ and for ${|z|<\bigl(np^{1/(p-1)}\bigr)^{-1}\sigma(R_P)}$ if ${p\le n}$. 
\end{theorem}
\paragraph{Proof}
For the case ${n < p\le \infty}$ see Schmidt \cite[Lemma~1]{Sc90}. As indicated by Schmidt,
the case ${n < p< \infty}$ is a direct consequence of a result of Dwork and Robba~\cite{DR79}.
The case ${p\le n}$  is due to Dwork and van der Poorten \cite[Theorem~3]{DP92}, who confirmed a conjecture of Schmidt. \qed

\medskip
Thus,  assuming~$P$ to be $w$-separable, everywhere above one can take ${\rho= c(p,n)^{-1}\sigma(R_P)}$, where
\begin{equation}
\label{ecpn}
c(p,n) =
\begin{cases}
1, & n<p\le \infty,\\
np, & p\le n.
\end{cases}
\end{equation}
Put
\begin{equation}
\label{esig}
\Sigma =
\begin{cases}
\max\bigl\{\sigma(R_P)^{-1}, (6|P|)^n\bigr\}, & p=\infty, \\
\max\bigl\{c(p,n)\sigma(R_P)^{-1}, |P|^n\bigr\}, & p<\infty,
\end{cases}
\end{equation}
The following theorem is now immediate.

\begin{theorem}[Local Eisenstein Theorem, regular case]
\label{tloreg}
Let ${f(z)\in \C_p[[z]]}$ written as in~\eqref{eserf} satisfy the polynomial equation ${P\bigl(z,f(z)\bigr)=0}$, where the polynomial ${P(z,w)\in\C_p[z,w]}$ is not divisible by~$z$ and $w$-separable. We normalize~$P$ so the polynomial ${P(0,w)\in \C_p[w]}$ is  monic.  
Then for ${k=0, 1, 2,\ldots}$ we have
${|a_k|\le 
3|P|\Sigma^k}$ when  ${p=\infty}$, and
${|a_k|\le |P|\Sigma^k}$ when ${p<\infty}$. \qed
\end{theorem}

For applications it is convenient to replace $\sigma(R_P)$, which is defined in terms of the roots of~$R_P$, by a quantity defined in terms of its coefficients. Let~$\mu$ be the order of~$0$ as the root of $R_P(z)$ and~$\gamma$ its lowest non-zero coefficient, so that   
\begin{equation}
\label{emugam}
R_P(z)=\gamma z^\mu+\text{higher powers of~$z$}. 
\end{equation}
Proposition~\ref{pboundroot} implies that
\begin{equation}
\label{esigvarsig}
\sigma(R_P)^{-1}\le 
\begin{cases}
2|R_P/\gamma|, & p=\infty,\\
 |R_P/\gamma|, & p<\infty.
 \end{cases}
\end{equation}
We obtain the following statement.

\begin{corollary}
\label{cloreg}
Let ${f(z)\in \C_p[[z]]}$ written as in~\eqref{eserf} satisfy the polynomial equation ${P\bigl(z,f(z)\bigr)=0}$, where the polynomial ${P(z,w)\in\C_p[z,w]}$ is not divisible by~$z$ and $w$-separable. We normalize~$P$ so the polynomial ${P(0,w)\in \C_p[w]}$ is  monic.  
Put
\begin{equation}
\label{exi}
\Xi =
\begin{cases}
\max\bigl\{2|R_P/\gamma|, (6|P|)^n\bigr\}, & p=\infty, \\
\max\bigl\{c(p,n)|R_P/\gamma|, |P|^n\bigr\}, & p<\infty,
\end{cases}
\end{equation}
Then for ${k=0, 1, 2,\ldots}$ we have
${|a_k|\le 
3|P|\Xi^k}$ when  ${p=\infty}$, and
${|a_k|\le |P|\Xi^k}$ when ${p<\infty}$. \qed
\end{corollary}

\section{Local Eisenstein Theorem: the General Case}
\label{slogen}
We now allow the series~$f$ to admit a finite pole at~$0$ and finite ramification:
\begin{equation}
\label{efram}
f(z) = \sum_{k=\kappa }^\infty a_kz^{k/e} \in \C_p((z^{1/e})),
\end{equation}
where~$e$ is a positive integer and~$\kappa $ is an integer which may be positive or negative or~$0$. At the moment we do not assume that~$\kappa $ is maximal possible and that~$e$ is minimal possible (that is, it may well happen that ${a_{\kappa }= 0}$ and/or that ${f\in \C_p((z^{1/e'}))}$ for some ${e'<e}$).  

We need to introduce one technical notion. Let $P(z,w)$ be a non-zero polynomial over some  field, and~$k$ an integer (which may be negative). Then there exists a unique integer~$N$ such 
$$
P_k(z,w):=z^NP(z,z^kw)
$$
is a polynomial not divisible by~$z$. The polynomial~$P$ will be called \textsl{$k$-normalized} if the polynomial $P_k(0,w)$ is monic. This is a \textsl{moderate normalization}, as defined in Section~\ref{sprelim}. The normalization hypothesis of Theorem~\ref{tloreg} means exactly that~$P$ is $0$-normalized.

It might be worth remarking that the polynomial~$P_k$ has the same set of non-zero coefficients as~$P$; in particular ${|P_k|=|P|}$. Also, the $w$-resultant of~$P_k$ is equal to the $w$-resultant of~$P$ times a power of~$z$. In particular, the quantity~$\Sigma$, defined in~\eqref{esig}, as well as the quantity~$\Sigma_e$ defined in~\eqref{esige} below, are the same for~$P$ and~$P_k$.

Now let~$f$ be as in the beginning of this section, and ${P(z,w)\in \C_p[z,w]}$  a $w$-separable polynomial of $w$-degree~$n$ such that ${P\bigl(z,f(z)\bigr)=0}$, and $R_P(z)$ be the $w$-resultant of $P(z,w)$ and $P'_w(z,w)$. Put
\begin{equation}
\label{esige}
\Sigma_e =
\begin{cases}
\max\bigl\{\sigma(R_P)^{-1}, (6|P|)^n\bigr\}, & p=\infty, \\
\max\bigl\{c(p,n)^e\sigma(R_P)^{-1}, |P|^n\bigr\}, & p<\infty,
\end{cases}
\end{equation}
where ${\sigma(\cdot)}$ is defined in~\eqref{esig}. 
If ${e=1}$ (no ramification) and~$P$ is $\kappa $-normalized, then we have the estimates
${|a_k|\le 
3|P|\Sigma_1^{k-\kappa }}$ when  ${p=\infty}$, and
${|a_k|\le |P|\Sigma_1^{k-\kappa }}$ when ${p<\infty}$. This follows by applying Theorem~\ref{tloreg} to the series ${z^{-\kappa }f(z)}$ and the polynomial $P_{\kappa }(z,w)$.

We believe that in the general case, after a suitable moderate normalization of~$P$, the estimates ${|a_k|\le 
3|P|\Sigma_e^{(k-\kappa )/e}}$ when  ${p=\infty}$, and
${|a_k|\le |P|\Sigma_e^{(k-\kappa )/e}}$ when ${p<\infty}$  must hold.
Unfortunately, we can prove only a slightly weaker result.

\begin{theorem}[Local Eisenstein Theorem, general case]
\label{tlogen}
Let $f(z)$ as in~\eqref{efram} satisfy the polynomial equation ${P\bigl(z,f(z)\bigr)=0}$, where the polynomial ${P(z,w)\in\C_p[z,w]}$ is $w$-separable and $\lfloor \kappa /e\rfloor$-normalized. Then for ${k\ge \kappa }$ we have the estimates ${|a_k|\le 
3|P|\Sigma_e^{k/e-\lfloor \kappa /e\rfloor}}$ when  ${p=\infty}$, and
${|a_k|\le |P|\Sigma_e^{k/e-\lfloor \kappa /e\rfloor}}$ when ${p<\infty}$.
\end{theorem}

Replacing $f(z)$ and $P(z,w)$ by ${z^{-\lfloor \kappa /e\rfloor}f(z)}$ and $P_{\lfloor \kappa /e\rfloor}(z,w)$ (as indicated above, this does affect the quantity~$\Sigma_e$), we may assume that ${\lfloor \kappa /e\rfloor=0}$; in particular, ${\kappa \ge 0}$. Defining ${a_k=0}$ for ${0\le k<\kappa }$, we reduce the theorem to the special case ${\kappa =0}$.
Thus,  we have to prove the following.

\begin{proposition}
In the set-up of Theorem assume that ${\kappa =0}$. Then for ${k\ge 0}$ we have the estimates ${|a_k|\le 
3|P|\Sigma_e^{k/e}}$ when  ${p=\infty}$, and
${|a_k|\le |P|\Sigma_e^{k/e}}$ when ${p<\infty}$.
\end{proposition}

\paragraph{Proof}
It is quite analogous to that of Theorem~\ref{tloreg}. Put ${\tilf(z)=f(z^e)}$ and ${\tilP(z,w)=P(z^e,w)}$, so that ${\tilP\bigl(z,\tilf(z)\bigr)=0}$. We have ${|\tilP|=|P|}$ and ${\sigma(R_\tilP)=\sigma(R_P)^{1/e}}$. Theorem~\ref{tdrsv} implies now that $\tilf(z)$ converges in the circle ${|z|<\tilrho:=c(p,n)^{-1}\sigma(R_P)^{1/e}}$. Put
$$
\tilvarrho =
\begin{cases}
\min\left\{\tilrho, (6|P|)^{-n/e}\right\}, & p=\infty,\\
\min\left\{\tilrho, |P|^{-n/e}\right\}, & p<\infty.
\end{cases}
$$
The polynomial ${q(w):=P(0,w)}$ is monic by the assumption, and, in the notation of~\eqref{equs}, we have
$$
\tilP(z,w)=q(w)+z^eQ(z^e,w).
$$
Arguing exactly as in the proof of Proposition~\ref{pbor}, we obtain that for ${|z|<\tilvarrho}$  we have ${|f(z)| \le 3|P|}$ if ${p=\infty}$, and ${|f(z)| \le |P|}$ if ${p<\infty}$. Now Corollary~\ref{ccoef} implies
$$
|a_k|\le 
\begin{cases}
3|P|\tilvarrho^{-k}, & p=\infty, \\
|P|\tilvarrho^{-k}, & p<\infty.
\end{cases}
$$
Since ${\tilvarrho^{-1}=\Sigma_e^{1/e}}$, this completes the proof.\qed

\medskip

We also state an analog of Corollary~\ref{cloreg}, with $\sigma(R_P)^{-1}$ replaced by $|R_P/\gamma|$. 

\begin{corollary}
\label{clogen}
Let $f(z)$ as in~\eqref{efram} satisfy the polynomial equation ${P\bigl(z,f(z)\bigr)=0}$, where the polynomial ${P(z,w)\in\C_p[z,w]}$ is $w$-separable and $\lfloor \kappa /e\rfloor$-normalized. 
Put
\begin{equation}
\label{exie}
\Xi_e =
\begin{cases}
\max\bigl\{2|R_P/\gamma|, (6|P|)^n\bigr\}, & p=\infty, \\
\max\bigl\{c(p,n)^e|R_P/\gamma|, |P|^n\bigr\}, & p<\infty,
\end{cases}
\end{equation}
Then for ${k\ge \kappa }$ we have the estimates ${|a_k|\le 
3|P|\Xi_e^{k/e-\lfloor \kappa /e\rfloor}}$ when  ${p=\infty}$, and
${|a_k|\le |P|\Xi_e^{k/e-\lfloor \kappa /e\rfloor}}$ when ${p<\infty}$. \qed 
\end{corollary}

\section{Global Eisenstein Theorem}
\label{sglo}
In this section~$K$ is a number field of degree ${d=[K:\Q]}$ and  $M_K$ is the set of its absolute values normalized as indicated in Section~\ref{sprelim}.

Recall that by an \textsl{$M_K$-divisor} we mean associating to every ${v\in M_K}$ a positive real number~$A_v$ with the following property: for all but finitely many~$v$ we have ${A_v=1}$. The $M_K$-divisor is \textit{effective} if ${A_v\ge 1}$ for all~$v$. The \textsl{height} of an $M_K$-divisor ${\calA=(A_v)_{v\in M_K}}$ is defined as 
$$
\vheight(\calA)= d^{-1}\sum_{v\in M_K} d_v\log^+A_v,
$$
where ${d_v=[K_v:\Q_v]}$ is the absolute local degree of~$v$. For an effective divisor $\log^+$ can be replaced by $\log$.

{\sloppy

\begin{theorem}[Global Eisenstein Theorem, regular case]
\label{tgloreg}
Let ${P(z,w)\in K[z,w]}$ be a $w$-separable polynomial of degrees 
${\deg_zP=m}$ and  ${\deg_wP=n}$, and
$$
f(z)=\sum_{k=0}^\infty a_kz^k \in \bar K[[z]]
$$
a power series  satisfying ${P\bigl(z,f(z)\bigr)=0}$. Then there exist effective $M_K$-divisors ${\calA'=(A_v')_{v\in M_K}}$ and ${\calA=(A_v)_{v\in M_K}}$ such that 
$$
|a_k|_v\le A_v'A_v^k\qquad (k=0,1,2,\ldots)
$$ 
for any ${v\in M_K}$  anyhow extended to $\bar K$,  and such that 
\begin{equation}
\label{eeisenreg}
\vheight(\calA') \le \pheight(P)+\log 3,  \quad  \vheight(\calA)\le (3n-1)\pheight(P) + 3n\log (mn) +7n. 
\end{equation}
 \end{theorem}
 
 }

\begin{remark}
The $w$-separability assumption can be dropped for the price of slightly increasing the estimates~\eqref{eeisenreg}. Indeed,  the classical  ``Gauss-Mahler-Gelfond'' lemma (see, for instance, \cite[Theorem~1.6.13]{BG06}) implies that if~$f$,~$g$ are polynomials in~$r$ variables ${x_1, \ldots, x_r}$ with algebraic coefficients, and ${g\mid f}$, then ${\pheight(g)\le \pheight(f)+n_1+\cdots+n_r}$, where ${n_i=\deg_{x_i}f}$. Hence, if the polynomial~$P$ is not $w$-separable, then we may replace it by its square free part, which is $w$-separable and whose height is at most ${\pheight(P)+m+n}$. 
\end{remark}

\paragraph{Proof}
We may assume that the polynomial~$P$ is not divisible by~$z$ and is normalized  in such a way  that the $w$-polynomial $P(0,w)$ is monic ($0$-normalized in terminology of Section~\ref{slogen}). For every ${v\in M_K}$ we define~$A_v'$ and~$A_v$ as Corollary~\ref{cloreg} suggests:
$$
A_v =
\begin{cases}
\max\bigl\{2|R_P/\gamma|_v, (6|P|_v)^n\bigr\}, & v\mid\infty, \\
\max\bigl\{c(p,n)|R_P/\gamma|_v, |P|_v^n\bigr\}, & v\mid p<\infty,
\end{cases}\qquad
A_v'=
\begin{cases}
3|P|_v, & v\mid\infty, \\
|P|_v, & v\mid p<\infty,
\end{cases}
$$
where $R_P(z)$ is the $w$-resultant of $P(z,w)$ and $P'_w(z,w)$, written as in~\eqref{emugam}, and $c(p,n)$ is defined in~\eqref{ecpn}.
Since the polynomial~$P$ has a coefficient equal to~$1$, both divisors~$\calA'$ and~$\calA$ are effective. Theorem~\ref{tloreg} implies that ${|a_k|_v\le A_v'A_v^k}$ however~$v$ is extended to~$\bar K$. 

Clearly ${\vheight(\calA') \le \pheight(P)+\log 3}$. Now, 
$$
\log A_v \le 
\begin{cases}
\log |R_P/\gamma|_v+ n\log |P|+n\log 6, & v\mid\infty, \\
\log |R_P/\gamma|_v+ n\log |P|+\log c(p,n), & v\mid p<\infty,
\end{cases}
$$
It follows that
\begin{equation}
\label{eintermi}
\vheight(\calA) \le \pheight(R_P) +\log 2 + n\pheight(P)+ n\log 6 + \sum_{p<\infty}\log c(p,n).
\end{equation}
For the latter sum we have 
$$
\sum_{p<\infty}\log c(p,n)= \pi(n)\log n+\sum_{p\le n}\log p  \le 2.3n,
$$
where we use inequalities~(3.6) and~(3.32) from \cite[pages 69--71]{RS62}. Combining all this with Proposition~\ref{prp}, we obtain, after a little calculation, the wanted estimate 
$$
\vheight(\calA)\le (3n-1)\pheight(P) + 3n\log (mn)  +7n.
\eqno\square
$$

\begin{theorem}[Global Eisenstein Theorem, general case]
\label{tglogen}
Let
$$
f(z)=\sum_{k=\kappa }^\infty a_kz^{k/e} \in \bar K((z^{1/e}))
$$
be a power series over~$K$ satisfying ${P\bigl(z,f(z)\bigr)=0}$.
Then there exist effective $M_K$-divisors ${\calA'=(A_v')_{v\in M_K}}$ and ${\calA=(A_v)_{v\in M_K}}$ such that ${|a_k|_v\le A_v'A_v^{k/e-\lfloor \kappa /e\rfloor}}$ for ${v\in M_K}$ anyhow extended to~$\bar K$ and  for ${k\ge \kappa }$, and 
\begin{equation}
\label{eeisengen}
\vheight(\calA') \le \pheight(P)+\log 3,  \quad  \vheight(\calA)\le (3n-1)\pheight(P) + 3n\log (mn) +7en. 
\end{equation}
\end{theorem}

\paragraph{Proof}
The proof is identical to that of Theorem~\ref{tgloreg}, but now instead of Corollary~\ref{cloreg} one should use Corollary~\ref{tlogen}; in particular, we define
\begin{equation}
\label{ea_v}
A_v =
\begin{cases}
\max\bigl\{2|R_P/\gamma|_v, (6|P|_v)^n\bigr\}, & v\mid\infty, \\
\max\bigl\{c(p,n)^e|R_P/\gamma|_v, |P|_v^n\bigr\}, & v\mid p<\infty,
\end{cases}\qquad
A_v'=
\begin{cases}
3|P|_v, & v\mid\infty, \\
|P|_v, & v\mid p<\infty,
\end{cases}
\end{equation} 
where~$P$ is assumed to be $\lfloor\kappa/e\rfloor$-normalized. 
We leave the further details to the reader.\qed

\medskip

We conclude this section with one application of Theorem~\ref{tglogen},  which will be used in~\cite{BSS12}. 
Eisenstein theorem implies that for all but finitely ${v\in M_K}$ the $v$-adic norms of all the coefficients of an algebraic power series are bounded by~$1$. We want to estimate the size of the  finite set of exceptional~$v$ for which this fails. 

Let us make some definitions. We define the \textsl{absolute norm} $\norm v$ of ${v\in M_K}$ as the absolute norm of the corresponding prime ideal if ${v\mid p<\infty}$; for ${v\mid \infty}$ we set ${\norm v=1}$. 

We define the \textsl{height} of a finite subset ${S\subset M_K}$ as 
${\vheight(S) = d^{-1}\sum_{v\in S} \log \norm v}$. 
(Recall that ${d=[K:\Q]}$.) 

So far we dealt with an individual series~$f$. However, if the polynomial ${P(z,w)\in K[z,w]}$ is $w$-separable, then the Puiseux theorem implies existence of ${n=\deg_wP}$ distinct series ${f_1, \ldots, f_n}$, which can be written as
\begin{equation}
\label{ef_i}
f_i(z)=\sum_{k=\kappa_i}^\infty a_{ik}z^{k/e_i} \qquad (i=1, \ldots, n),
\end{equation}
and which satisfy ${P\bigl(z,f_i(z)\bigl)=0}$. 

{\sloppy

\begin{theorem}
Let
${P(z,w)\in K[z,w]}$ be a $w$-separable polynomial with ${m=\deg_zP}$ and ${n=\deg_wP}$, and let ${f_1, \ldots, f_n}$ be the~$n$ distinct series, written as in~\eqref{ef_i} and satisfying ${P\bigl(z,f_i(z)\bigl)=0}$.  Let~$S$ be the (finite) set of ${v\in M_K}$ such that ${|a_{ik}|_v>1}$ for some coefficient $a_{ik}$ and some extension of~$v$ to~$\bar K$. Then
\begin{equation}
\label{ehes}
\vheight(S) \le 3n\bigl(\pheight(P)+\log(mn)+1). 
\end{equation}
\end{theorem}

}

\paragraph{Proof}
For a non-archimedean ${v\in M_K}$  let~$\pi_v$ be a primitive element of~$v$ (a generator of the maximal ideal of the local ring of~$v$). Any ${\alpha \in K^\times}$ can be written as ${\alpha=\pi_v^\ell\eta}$ with ${\ell \in \Z}$ and~$\eta$ a $v$-adic unit. One verifies immediately that 
${\ell = d_v\log|\alpha|_v/\log \norm v}$ (where ${d_v=[K_v:\Q_v]}$ is the local degree),  which shows that for ${\alpha \in K^\times}$ the quotient ${d_v\log|\alpha|_v/\log \norm v}$ is an integer. In particular, if ${|\alpha|_v>1}$ then ${d_v\log|\alpha|_v\ge \log \norm v}$. 

Denote by $P_i(z,w)$ the $\lfloor\kappa_i/e_i\rfloor$-normalization of~$P$. As follows from Theorem~\ref{tglogen} together with definitions~\eqref{ea_v}, we may have ${|a_{ik}|_v>1}$ (for some extension of~$v$ to~$K$) only if either ${|P_i|_v>1}$, or ${|R_P/\gamma|_v>1}$, or ${v\mid p\le n}$ or ${v\mid \infty}$. Since each of the numbers ${|P_i|_v}$, ${|R_P/\gamma|_v}$ and ${|p^{-1}|_v}$ is equal to ${|\alpha|_v}$ for some ${\alpha\in K}$, we obtain 
$$
\log\norm v \le 
\begin{cases}
d_v\log |P_i|_v, &|P_i|_v>1,\\
d_v\log |R_P/\gamma|_v, &|R_P/\gamma|_v>1,\\
d_v\log p, & v\mid p. 
\end{cases}
$$
It follows that
\begin{align*}
\vheight(S) &\le \pheight(R_P/\gamma)+ \pheight(P_1) +\cdots+ \pheight(P_n) +\sum_{p\le n} \log p\\
& = \pheight(R_P)+ n\pheight(P) +\sum_{p\le n} \log p\\ 
&\le 3n\bigl(\pheight(P)+\log(mn)+1),
\end{align*}
where we use Proposition~\ref{prp} and estimate~(3.32) in \cite[page~71]{RS62}. \qed

\section{Estimates Involving the Initial Coefficient}
\label{sa0}
For certain types of applications one needs a result slightly subtler than Theorems~\ref{tgloreg}. Precisely, we want to express the divisor~$\calA'$ not in terms of the polynomial~$P$, but in terms of the initial coefficient~$a_0$. We prove the following theorem.

\begin{theorem}
\label{ta0reg}
In the set-up of Theorem~\ref{tgloreg}
there exists an effective $M_K$-divisor  ${\calA=(A_v)_{v\in M_K}}$ such that 
$$
|a_k|_v\le \max\{1,|a_0|_v\}A_v^k\qquad (k=0,1,2,\ldots)
$$ 
for any ${v\in M_K}$  anyhow extended to $\bar K$,  and such that 
\begin{equation}
\label{ea0reg}
 \vheight(\calA)\le 3n\pheight(P) + 3n\log (mn) +10n. 
\end{equation}
\end{theorem}

For the proof, we shall use the following modification of Corollary~\ref{ccoef}. 
\begin{proposition}
\label{pcoef}
In the set-up of Proposition~\ref{pbor}, the coefficients of the series~\eqref{eserf} satisfy
$$
|a_k|\le 
\begin{cases}
\max\{1,|a_0|\}(8|P|\varrho^{-1})^k, & p=\infty, \\
\max\{1,|a_0|\}(|P|\varrho^{-1})^k, & p<\infty
\end{cases}
\qquad (k=0,1,2,\ldots).
$$
\end{proposition}

\paragraph{Proof}
Using Corollary~\ref{ccoef}, we find that for ${|z|<\varrho}$  
$$
|f(z)-a_0|\le 
\begin{cases}
\frac{3|P||z|}{\varrho-|z|},& p=\infty, \\
|P||z|/\varrho, & p<\infty.
\end{cases}
$$
It follows that  in the archimedean case for ${|z|\le (1/2)\varrho |P|^{-1}}$ we have ${|f(z)-a_0|\le 3}$. Hence 
$$
|f(z)|\le 4\max\{1,|a_0|\}, 
$$
and applying~\eqref{eakle} with ${r=(1/2)\varrho |P|^{-1}}$, we obtain
${|a_k| \le 4\max\{1,|a_0|\}(2|P|\varrho^{-1} )^k}$, which implies  that ${|a_k|\le  \max\{1,|a_0|\}(8|P|\varrho^{-1})^k}$. 
Similarly, in the non-archimedean case for ${|z|<\varrho |P|^{-1}}$ we have ${|f(z)-a_0|<1}$. Hence 
${|f(z)|\le \max\{1,|a_0|\}}$,
and applying~\eqref{eakle} with arbitrary ${r<\varrho |P|^{-1}}$, we obtain
${|a_k| \le \max\{1,|a_0|\}(|P|\varrho^{-1} )^k}$.  \qed

\medskip

Now we have the following analog of Corollary~\ref{cloreg}.

\begin{corollary}
\label{ca0reg}
In the set-up of Corollary~\ref{cloreg} put 
$$
\Sigma =
\begin{cases}
8|P|\max\bigl\{2|R_P/\gamma|, (6|P|)^n\bigr\}, & p=\infty, \\
|P|\max\bigl\{c(p,n)|R_P/\gamma|, |P|^n\bigr\}, & p<\infty. 
\end{cases}
$$
Then ${|a_k|\le \max\{1,|a_0|\}\Sigma^k}$ for ${k=0,1,2,\ldots}$ \qed
\end{corollary}

\paragraph{Proof of Theorem~\ref{ta0reg}}
Same as the proof of Theorem~\ref{tgloreg}, but now, as Corollary~\ref{ca0reg} suggests, we put
$$
A_v =
\begin{cases}
8|P|_v\max\bigl\{2|R_P/\gamma|_v, (6|P|_v)^n\bigr\}, & v\mid\infty, \\
|P|_v\max\bigl\{c(p,n)|R_P/\gamma|_v, |P|_v^n\bigr\}, & v\mid p<\infty. 
\end{cases}
$$
Instead of~\eqref{eintermi} we have
\begin{equation}
\label{eintermia0}
\vheight(\calA) \le \pheight(R_P) +\log 2 + (n+1)\pheight(P)+ n\log 48 + \sum_{p<\infty}\log c(p,n).
\end{equation}
Arguing as in in the end of the proof of Theorem~\ref{tgloreg}, we find that the right-hand side of~\eqref{eintermia0} does not exceed ${3n\pheight(P) + 3n\log (mn) +10n}$.  \qed

\medskip

Similar results hold true in the general case as well. Here are the analogues of Corollary~\ref{clogen} and Theorem~\ref{tglogen}. We omit the details which are routine.

\begin{corollary}
In the set-up of Corollary~\ref{clogen}
put
\begin{equation*}
\Xi_e =
\begin{cases}
(8|P|)^e\max\bigl\{2|R_P/\gamma|, (6|P|)^n\bigr\}, & p=\infty, \\
|P|^e\max\bigl\{c(p,n)^e|R_P/\gamma|, |P|^n\bigr\}, & p<\infty,
\end{cases}
\end{equation*}
Then for ${k\ge \kappa }$ we have the estimate ${|a_k|\le 
\max\{1,|a_{e\lfloor\kappa/e\rfloor|}\}\Xi_e^{k/e-\lfloor \kappa /e\rfloor}}$. \qed 
\end{corollary}

Here we tacitly define ${a_k=0}$ for ${k<\kappa}$. In particular, ${\max\{1,|a_{e\lfloor\kappa/e\rfloor}|\}=1}$ if ${e\nmid\kappa}$.

\begin{theorem}
\label{ta0gen}
In the set-up of Theorem~\ref{tglogen}, there exists an effective $M_K$-divisor  ${\calA=(A_v)_{v\in M_K}}$ such that ${|a_k|_v\le \max\{1,|a_{e\lfloor\kappa/e\rfloor}|_v\}A_v^{k/e-\lfloor \kappa /e\rfloor}}$ for ${v\in M_K}$ anyhow extended to~$\bar K$ and any ${k\ge \kappa }$, and 
\begin{equation}
\label{ea0gen}
 \vheight(\calA)\le (3n+e-1)\pheight(P) + 3n\log (mn) +10en. 
\end{equation}
$ $\qed
\end{theorem}


\section{Fields Generated by the Coefficients}
\label{sfield}
If the polynomial $P(z,w)$ has coefficients in a field~$K$ (of characteristic~$0$), then any power series $f(z)$ satisfying ${P\bigl(z,f(z)\bigr)=0}$ has coefficients in a finite extension~$L$ of~$K$ of degree at most ${n=\deg_wP}$; this follows from the fact that there can be at most~$n$ distinct series $g(z)$ satisfying ${P(z,g(z))=0}$, and they include all the series obtained from~$f$ by the Galois conjugation over~$K$.

For applications  one needs to estimate the discriminant or some other invariants of the field~$L$ in the case when~$K$ is a number field;  see, for instance,~\cite{BC70,Bi97,Co70,Sc92}, where such estimates are crucial for Baker's method. 
In fact, Schmidt's interest in Eisenstein theorem was largely motivated by applications in Diophantine analysis~\cite{Sc92}, in particular, through estimating in~\cite{Sc91} the number fields generated by coefficients of certain algebraic power series.

The standard approach used in the articles quoted above   was to estimate the number of the coefficients of the series~$f$ needed to generate~$L$, and then to estimate the field~$L$ itself, in the form of estimating one of its generators, as in \cite{BC70,Sc91} or its discriminant, as in \cite[Lemma 2.4.2]{Bi97}. 

Here we follow similar approach, but introduce one technical novelty (see Lemma~\ref{lsum} below) which allows us to obtain results looking best possible for the method up to a constant factor.  Surprisingly, it turns out to be  more efficient to generate the field~$L$ coefficient by coefficient, passing through the subfields, rather than by all the coefficients  at once, as in \cite[Lemma 2.4.2]{Bi97}. In particular, we use item~\ref{isil} of Proposition~\ref{pdis} below only in the case when~$\bfa$ is a singleton.

Let us introduce some notation. Given an extension ${L/K}$ of number fields, we denote by $\partial_{L/K}$ the \textsl{normalized logarithmic relative discriminant}:
$$
\partial_{L/K}=\frac{\log\norm_{K/\Q}\calD_{L/K}}{[L:\Q]},
$$
where $\calD_{L/K}$ is the usual relative discriminant and $\norm_{K/\Q}$ is the norm map. The properties of this quantity are summarized in the following proposition, which will be used in the sequel without special reference. 
\begin{proposition}
\label{pdis}
\begin{enumerate}
\item
\label{iadd}
(additivity in towers)\quad  If ${K\subset L\subset M}$ is a tower of number fields, then  
${\partial_{M/K}= \partial_{L/K}+\partial_{M/L}}$.

\item
\label{ibase} 
(base extension)\quad
If~$K'$ is a finite extension of~$K$ and ${L'=LK'}$ then
${\partial_{L'/K'}\le \partial_{L/K}}$.

\item
\label{itri}
(triangle inequality)\quad If~$L_1$ and~$L_2$ are two extensions of~$K$, then 
${\partial_{L_1L_2/K}\le \partial_{L_1/K}+  \partial_{L_2/K}}$.

\item
\label{isil} 
(bounding in terms of the generators)\quad
Let ${\bfa=(a_1,\ldots,a_k)}$ be a point in ${\bar K^k}$.  Put ${L=K(\bfa)}$ and ${\nu=[L:K]}$. Then
$$
\partial_{L/K} \leq 2(\nu-1)\aheight(\bfa)+\log \nu. 
$$
\end{enumerate}
\end{proposition}

\paragraph{Proof}
Items~\ref{iadd} and~\ref{ibase} follow from the definition of the discriminant as the norm of the different, and the multiplicativity of the different in towers. Item~\ref{itri} is a direct consequence of the previous two. Item~\ref{isil} is due to Silverman \cite[Theorem 2]{Si84}. \qed

\medskip

In the sequel~$K$ is a number field, 
$$
P(z,w)= p_n(z) w^n+p_{n-1}(z) w^{n-1}+ \cdots + p_0(z) \in K[z,w]
$$ 
is a $w$-separable polynomial  with
$$
m=\deg_zP, \qquad n=\deg_wP,
$$
and ${f(z)\in \bar K((z^{1/e}))}$ is a power series satisfying ${P\bigl(z,f(z)\bigr)=0}$.

{\sloppy

More generally, since the polynomial $P(z,w)$ is $w$-separable,  there is~$n$ distinct series ${f_1=f, f_2, \ldots, f_n}$ satisfying ${P\bigl(z,f_i(z)\bigr)=0}$, with ${f_i(z) \in \bar K((z^{1/e_i}))}$ for some natural~$e_i$. 

}

We denote by~$L$ the number field generated over~$K$ be the coefficients of~$f$; more generally, we denote by ${L_1, \ldots, L_n}$ the number fields generated by the coefficients of ${f_1,  \ldots, f_n}$,  respectively.

\subsection{Integral case}

In this subsection  we consider the \textsl{integral case}, that is, we assume that  
$$
p_n(0)\ne 0. 
$$
This latter condition is equivalent to saying that the series ${f_1, \ldots, f_n}$ have no negative part: ${f_i(z)\in \bar K [[z^{1/e_i}]]}$. In particular, for ${f=f_1}$ we have 
$$
f(z)=\sum_{k=0}^\infty a_kz^{k/e}\in \bar K[[z^{1/e}]]. 
$$ 
Our main tool is Lemma~\ref{lsum} below. To state it, we need some more definitions. We denote by $\Lambda_k$ the field generated over~$K$ by the first~$k$ coefficients of~$f$; precisely,
$$
\Lambda_0=K, \qquad \Lambda_k= K(a_0, \ldots a_{k-1}) \qquad (k=1,2,3,\ldots). 
$$
Clearly, ${\Lambda_k=L}$ for sufficiently large~$k$. 
Further, put
${\lambda_k=[L:\Lambda_k]}$,
so that all but finitely many of~$\lambda_k$ are~$1$.

\begin{lemma}
\label{lsum}
In the set-up above, 
$$
\sum_{k=0}^\infty \frac ke(\lambda_k-\lambda_{k+1}) \le \ord_zP'_w\bigl(z,f(z)\bigr).
$$
\end{lemma}

\paragraph{Proof}
We may assume that ${f=f_1, f_2, \ldots, f_\nu}$ are the series obtained from~$f$ by Galois conjugation over~$K$. By the definition of the degrees~$\lambda_k$,  there is exactly ${\lambda_k-\lambda_{k+1}}$ indices ${i\in \{2, \ldots, \nu\}}$ satisfying ${\ord_z(f-f_i)=k/e}$. Hence
\begin{equation}
\label{euncomplete}
\sum_{k=0}^\infty \frac ke(\lambda_k-\lambda_{k+1}) = \ord_z\prod_{i=2}^\nu \bigl(f(z)-f_i(z)\bigr).
\end{equation}
Since the series ${f_1, \ldots, f_n}$ have no negative part, the product in the right-hand side of~\eqref{euncomplete} divides
$$
P'_w\bigl(z,f(z)\bigr) = p_n(z)\prod_{i=1}^n \bigl(f(z)-f_i(z)\bigr). 
$$
This completes the proof. 
 \qed

\medskip

Now we may state and prove the principal result of this section in the integral case.
Let ${D(z)=D_P(z)}$ be the $w$-discriminant of $P(z,w)$; it is not identically zero because~$P$ is $w$-separable. 

\begin{theorem}
\label{tintfield}
Assume that ${p_n(0)\ne 0}$. 
\begin{enumerate}
\item
\label{ione}
The field~$L$, generated by the coefficients of~$f$, satisfies 
$$
\partial_{L/K} \le 2(\nu-1)\aheight(a_0) +(8n-1)\,\ord_zP'_w\bigl(z,f(z)\bigr)\bigl( \pheight(P) + \log (mn)+3e \bigr),
$$ 
where ${\nu=[L:K]}$. 
\item
\label{iall}
Put ${E=\max\{e_1,\ldots,e_n\}}$. Then the number fields ${L_1, \ldots, L_n}$, generated over~$K$ by the coefficients of ${f_1, \ldots, f_n}$, respectively, satisfy
\begin{equation}
\label{esumord}
\sum_{i=1}^n\partial_{L_i/K} \le 2(n-1)\bigl(\pheight(P)+\log (n+1)\bigr) +(8n-1)\, \ord_zD(z)\bigl( \pheight(P) + \log (mn)+3E \bigr).
\end{equation}
\end{enumerate}
\end{theorem}

\paragraph{Proof}
Write 
${\mu_k=[\Lambda_{k+1}:\Lambda_k]=  \lambda_k/\lambda_{k+1}}$. 
Items~\ref{iadd} and~\ref{isil} of Proposition~\ref{pdis} imply that 
\begin{equation}
\label{eparpar}
\partial_{L/K} \le \sum_{k=0}^\infty\bigl( 2(\mu_k-1)\aheight(a_k) +\log \mu_k\bigr). 
\end{equation}
Theorem~\ref{ta0gen} gives an $M_L$-divisors~$\calA$, satisfying~\eqref{ea0gen}, and such that  
${\aheight(a_k) \le \aheight(a_0) + \vheight(\calA)k/e}$. 
Substituting this to~\eqref{eparpar}, we obtain
$$
\partial_{L/K} \le 2\aheight(a_0)\sum_{k=0}^\infty (\mu_k-1) + 2\vheight(\calA)\sum_{k=0}^\infty \frac ke(\mu_k-1) + \log \nu
$$
Lemma~\ref{lsum} implies that 
$$
\sum_{k=0}^\infty \frac ke(\mu_k-1)= \sum_{k=0}^\infty \frac ke\frac{\lambda_k-\lambda_{k+1}}{\lambda_{k+1}} \le \sum_{k=0}^\infty \frac ke(\lambda_k-\lambda_{k+1}) \le \ord_zP'_w\bigl(z,f(z)\bigr),
$$
Similarly,
$$
\sum_{k=0}^\infty (\mu_k-1)=  \sum_{k=0}^\infty \frac{\lambda_k-\lambda_{k+1}}{\lambda_{k+1}} \le \sum_{k=0}^\infty (\lambda_k-\lambda_{k+1})= \nu-1. 
$$ 
We obtain 
$$
\partial_{L/K} \le 2(\nu-1)\aheight(a_0) + 2\,\ord_zP'_w\bigl(z,f(z)\bigr) \vheight(\calA) + \log \nu. 
$$
Combining this with~\eqref{ea0gen}, item~\ref{ione} follows after a simplification.

Analogous inequalities hold  for every field~$L_i$. 
Since ${\nu_i=[L_i:K] \le n}$ and  
$$
D(z)= \prod_{i=1}^n P'_w\bigl(z,f_i(z)\bigr),
$$ 
summing up these~$n$ inequalities, we obtain
$$
\sum_{i=1}^n\partial_{L_i/K} \le 2(n-1)\bigl(\aheight(a_{1,0})+\cdots+\aheight(a_{n,0})\bigr) +(8n-1) \, \ord_zD(z)\bigl( \pheight(P) + \log (mn)+3E \bigr),
$$
where~$a_{i,0}$ is the initial coefficient of~$f_i$. It remains to notice that ${a_{1,0},\ldots,a_{n,0}}$ are the roots of the polynomial ${q(w)=P(0,w)}$, and item~\ref{iallroots} of Proposition~\ref{pmahler} gives
$$
\aheight(a_{1,0})+\cdots+\aheight(a_{n,0}) \le \pheight(q)+\log(n+1)\le \pheight(P)+\log(n+1). 
$$
This proves item~\ref{iall}. 
\qed

\medskip

Estimates in Theorem~\ref{tintfield} involve orders of vanishing, which is practical for certain applications, but not convenient to be used directly.  Therefore we give below  a ``pr\^et \`a porter'' version, ready to be used; it also gives a good idea of the quality of our estimates.   
\begin{corollary}
\label{cfriend}
Assume that ${p_n(0)\ne 0}$. 
\begin{enumerate}
\item
We have
\begin{equation}
\label{esumhe}
\sum_{i=1}^n\partial_{L_i/K} \le 16mn(n-1)\bigl( \pheight(P) + \log (mn)+3E\bigr). 
\end{equation}

\item
Assume that the field ${L=L_1}$ is of degree~$\nu$ over~$K$. Then
\begin{equation}
\label{eindiv}
\partial_{L/K} \le \frac1\nu16mn(n-1)\bigl( \pheight(P) + \log (mn)+3E\bigr). 
\end{equation}

\end{enumerate}

\end{corollary}

\paragraph{Proof}
Since ${\ord_zD(z)\le \deg D(z)\le 2m(n-1)}$, estimate~\eqref{esumhe} follows from~\eqref{esumord} after easy transformations.  Next, if ${[L:K]=\nu}$ then among the fields~$L_i$ there are~$\nu$ fields conjugate to~$L$ over~$K$. If~$L_i$ is conjugate to~$L$ then ${\partial_{L_i/K}=\partial_{L/K}}$. Hence the left-hand side of~\eqref{esumhe} is at least $\nu\partial_{L/K}$, which proves~\eqref{eindiv}.\qed

\medskip

We believe that the order of magnitude in the estimates~\eqref{esumhe} and~\eqref{eindiv} is best possible, but the numerical constant~$16$ can probably be replaced by~$8$.

\subsection{General case}
In this subsection we consider the general case; that is, we no longer assume that ${p_n(0)\ne 0}$. One can treat it similarly, using the general version of the Eisenstein theorem. But it turns out to be more practical to reduce it to the integral case treated above.

\begin{theorem}
The number fields ${L_1, \ldots, L_n}$, generated over~$K$ by the coefficients of ${f_1, \ldots, f_n}$, respectively, satisfy
\begin{equation}
\label{esumordgen}
\sum_{i=1}^n\partial_{L_i/K} \le 2(n-1)(\pheight(P)+4n) +(8n-1)\,\ord_zD(z)\bigl( \pheight(P) + 5n +\log m\bigr).
\end{equation}
\end{theorem}

\paragraph{Proof}
We may assume that the polynomial $P(z,w)$ is not divisible by~$z$.  
There exists an algebraic number~$\zeta$, which is either~$0$ or an $n$-th root of unity,  and which is distinct from any root of the polynomial $P(0,w)$.  Writing the polynomial
${\Q(z,w)= w^nP(z,w^{-1}+\zeta)}$ as 
$$
Q(z,w)=q_n(z)w^n+ q_{n-1}(z)w^{n-1}+\cdots+q_0(z),
$$ 
we find ${q_n(z)=P(z,\zeta)}$, and, in particular, ${q_n(0)\ne0}$. 
The series  ${g_i(z)=(f_i(z)-\zeta)^{-1}}$ satisfies ${Q\bigl(z,g_i(z)\bigr)=0}$, and its coefficients generate the field $L_i(\zeta)$ over $K(\zeta)$. We may also notice that the polynomials~$P$ and~$Q$ have the same $w$-discriminant (up to the sign); as before, we denote this discriminant by $D(z)$. 

We may now apply item~\ref{iall} of Theorem~\ref{tintfield} to the polynomial $Q(z,w)$, series $g_i(z)$ and the fields $L_i(\zeta)/K(\zeta)$. We obtain
\begin{equation*}
\sum_{i=1}^n\partial_{L_i(\zeta)/K(\zeta)} \le 2(n-1)\bigl(\pheight(Q)+\log (n+1)\bigr) +(8n-1) \, \ord_zD(z)\bigl( \pheight(Q) + \log (mn)+3E \bigr). 
\end{equation*}
By the choice of~$\zeta$ we have ${\aheight(\zeta)=0}$ and ${[K(\zeta):K]\le n-1}$. Proposition~\ref{ptrah} implies that \begin{equation*}
\pheight(Q)\le \pheight(P)+n\log 2+\log(n+1), 
\end{equation*}
and Proposition~\ref{pdis} implies that 
\begin{equation*}
\partial_{L_i/K}\le \partial_{L_i(\zeta)/K(\zeta)}+\partial_{K(\zeta)/K} \le \partial_{L_i(\zeta)/K(\zeta)}+\log (n-1). 
\end{equation*}
Combining the last three inequalities, we obtain~\eqref{esumordgen} after an obvious transformation. 
\qed

\medskip

We again give a ``pr\^et \`a porter'' version; the proof is the same as for Corollary~\ref{cfriend} and is left out.

\begin{corollary}
Assume that ${p_n(0)\ne 0}$. 
\begin{enumerate}
\item
We have
\begin{equation*}
\sum_{i=1}^n\partial_{L_i/K} \le 16mn(n-1)\bigl( \pheight(P) + 5n +\log m\bigr). 
\end{equation*}

\item
Assume that the field ${L=L_1}$ is of degree~$\nu$ over~$K$. Then
$$
\partial_{L/K} \le \frac2\nu16mn(n-1)\bigl( \pheight(P) + 5n + \log m\bigr). 
\eqno\square
$$
\end{enumerate}
\end{corollary}

\bigskip

\noindent
\textbf{Yuri Bilu}\\
IMB, Universit\'e Bordeaux 1\\
351 cours de la Libération\\
33405 Talence CEDEX\\
France

\bigskip

\noindent
\textbf{Alexander Borichev}\\
CMI, Universit\'e de Provence\\
Technop\^ole Ch\^ateau-Gombert\\
39, rue F. Joliot Curie\\
13453 Marseille Cedex 13\\
France 

\bigskip

\end{document}